\documentclass[twoside]{article}
\usepackage{amsmath,amsthm,amssymb}
\usepackage{url}
\usepackage[spanish,english]{babel}
\usepackage{fancyhdr}
\fancypagestyle{plain}{\fancyhf{}%

}

\fancypagestyle{headings}{\fancyhf{}%
\fancyhead[LE,RO]{\thepage}%
\fancyhead[RE]{Jos\'e Mijares, Jes\'us Nieto}%
\fancyhead[LO]{A parametrization of the abstract Ramsey theorem}%
\fancyfoot[C]{\footnotesize{ \pageref{begin-art}--\pageref{end-art}}}}

\newtheorem{clm}{Claim}
\newtheorem{thm}{Theorem}
\newtheorem{lem}{Lemma}
\newtheorem{coro}{Corollary}
\theoremstyle{definition}
\newtheorem{defn}{Definition}
\newtheorem{prop}{Proposition}

\begin{document}
\label{begin-art}
\pagestyle{headings}
\thispagestyle{plain} 
\selectlanguage{english}
\begin{center}
{\LARGE\bfseries A parametrization of the abstract Ramsey theorem \par }

\vspace{5mm}

{\large Jos\'e G. Mijares (\url{jmijares@euler.ciens.ucv.ve})}\\[1mm]
Departamento de Matem\'aticas IVIC\\ Escuela de Matem\'atica\\
Universidad Central de Venezuela\\
{\large Jes\'us E. Nieto (\url{jnieto@usb.ve})}\\[1mm]
Departamento de Matem\'aticas\\
Universidad Sim\'on Bol\'\i var
\end{center}

\begin{abstract}
We give a parametrization with perfect subsets of $2^{\infty}$ of the abstract Ramsey theorem (see \cite{todo}) Our main tool is an extension of the parametrized version of the combinatorial forcing developed in \cite{nash} and \cite{todo}, used in \cite{mij} to the obtain a parametrization of the abstract Ellentuck theorem. As one of the consequences, we obtain a parametrized version of the Hales-Jewett theorem. Finally, we conclude that the family of perfectly ${\cal S}$-Ramsey subsets of $2^{\infty}\times {\cal R}$ is closed under the Souslin operation.
\end{abstract}

\selectlanguage{english}
\section{Introduction}

\noindent In \cite{todo}, S. Todorcevic presents an abstract characterization of those topological spaces in which an analog of Ellentuck's theorem \cite{ellen} can be proven. These are called \emph{topological Ramsey spaces} and the main result about them is referred to in \cite{todo} as \emph{abstract Ellentuck theorem}. In \cite{mij}, a parametrization with perfect subsets of $2^{\infty}$ of the abstract Ellentuck theorem is given, obtaining in this way new proofs of parametrized versions of the Galvin-Prikry theorem \cite{galpri} (see \cite{miltodo}) and of Ellentuck's theorem (see \cite{pawl}), as well as a parametrized version of Milliken's theorem \cite{mill}. The methods used in \cite{mij} are inspired by those used in \cite{farah} to obtain a parametrization of the semiselecive version of Ellenctuck's theorem. 

\medskip

\noindent Nevertheless, topological Ramsey spaces are a particular kind of a more general type of spaces (introduced in \cite{todo}), in which the Ramsey property can be characterized in terms of the abstract Baire property. These are called \emph{Ramsey spaces}. One of such spaces, known as the \emph{Hales-Jewett space}, is described bellow (for a more complete description of this -- non topological-- Ramsey space, see \cite{todo}). S. Todorcevic has given a characterization of Ramsey spaces which is summed up in a result known as the \emph{abstract Ramsey theorem}. It tunrs out that the abstract Ellentuck theorem is a consequence of the abstract Ramsey theorem (see \cite{todo}). Definitions of all these concepts will be given bellow.

\medskip

\noindent In this work we adapt in a natural way the methods used in \cite{mij} in order to obtain a parametrized version of the abstract Ramsey theorem. In this way, we not only generalize the results obtained in \cite{mij} but we also obtain, in corollary \ref{paraHJ} bellow, a parametrization of the infinite dimensional version of the Hales-Jewett theorem \cite{HaJ} (see \cite{todo}), which is the analog to Ellentuck's theorem corresponding to the Hales-Jewett space. 

\medskip

\noindent In the next section we summarize the definitions and main results related to Ramsey spaces given in \cite{todo}. In section \ref{paramet} we introduce the (parametrized) combinatorial forcing adapted to the context of Ramsey spaces and present our main result (theorem \ref{main} bellow). Finally, we conclude that the generalization of the \emph{perfectly Ramsey property} (see \cite{dip} and \cite{pawl}) to the context of Ramsey spaces is preserved by the Souslin operation (see corollary \ref{souslin} bellow).

\medskip

\noindent We'll use the following definitions and results concerning to perfect sets and trees (see \cite{pawl}). For $x=(x_n)_n\in 2^{\infty}$, $x_{|k}=(x_0,x_1,\dots,x_{k-1})$. For $u\in 2^{\<infty}$, let $[u]=\{x\in 2^{\infty}\colon (\exists k)(u=x_{|k})\}$ and denote the \emph{length} of $u$ by $|u|$. If $Q\subseteq2^{\infty}$ is a perfect set, we denote $T_Q$ its asociated perfect tree. For $u$, $v=(v_0,\dots,v_{|v|-1})\in 2^{<\infty}$, we write $u\sqsubseteq v$ to mean $(\exists k\leq |v|)(u=(v_0,v_1,\dots,v_{k-1}))$. Given $u\in 2^{<\infty}$, let $Q(u)=Q\cap [u(Q)]$, where $u(Q)$ is defined as follows: $\emptyset(Q)=\emptyset$. If $u(Q)$ is already defined, find $\sigma \in T_Q$ such that $\sigma$ is the $\sqsubseteq$-extension of $u(Q)$ where the first ramification occurs. Then, set $(u^{\smallfrown}i)(Q)=\sigma^{\smallfrown}i$, $i=1,0$. Where "$^{\smallfrown}$" is concatenation. Thus, for each $n$, $Q=\bigcup\{Q(u)\colon u\in2^n\}$. For $n\in\mathbb{N}$ and perfect sets $S$, $Q$, we write $S\subseteq_n Q$ to mean $S(u)\subseteq Q(u)$ for every $u\in 2^{\infty}$. Thus "$\subseteq_n$" is a partial order and, if we have chosen $S_u\subseteq Q(u)$ for every $u\in 2^n$, then $S=\bigcup_u S_u$ is perfect, $S(u)=S_u$ and $S\subseteq_n Q$. The \emph{property of fusion} of this order is: if $Q_{n+1}\subseteq_{n+1}Q_n$ for $n\in \mathbb{N}$, then $Q=\cap_nQ_n$ is perfect and $Q\subseteq_n Q_n$ for each $n$.    

\section{Abstract Ramsey theory}
\noindent We introduce some definitions and results due to Todorcevic (see \cite{todo}). Our objects will be structures of the form $({\cal R},{\cal S},\leq,\leq^0,r,s)$ where $\leq$ and $\leq^0$ are relations on ${\cal S}\times {\cal S}$ and ${\cal R}\times{\cal S}$ respectively; and $r$, $s$ give finite approximations:
$$r\colon {\cal R}\times \omega \to {\cal AR}\qquad s\colon {\cal S}\times \omega \to {\cal AS}$$we denote $r_n(A)=r(A,n)$, $s_n(X)=s(X,n)$, for $A\in {\cal R}$, $X\in {\cal S}$, $n\in \mathbb{N}$. The following three axioms are assumed for every $({\cal P},p)\in \{({\cal R},r),({\cal S},s)\}$.

\begin{itemize}
\item[(A.1)] $p_0(P)=p_0(Q)$, for all $P$, $Q\in {\cal P}$. 
\item[(A.2)] $P\neq Q\Rightarrow p_n(P)\neq p_n(Q)$ for some $n\in \mathbb{N}$.
\item[(A.3)] $p_n(P)=p_m(Q)\Rightarrow n=m$ and $p_k(P)=p_k(Q)$ if $k<n$.
\end{itemize}   

\noindent In this way we can consider elements of ${\cal R}$ and ${\cal S}$ as infinite sequences $(r_n(A))_{n\in \mathbb{N}}$, $(s_n(X))_{n\in \mathbb{N}}$. Also, if $a\in {\cal AR}$ and $x\in{\cal AS}$ we can think of $a$ and $x$ as finite sequences $(r_k(A))_{k<n}$, $(s_k(X))_{k<m}$ respectively; with $n$, $m$ the unique integers such that $r_n(A)=a$ and $s_m(X)=x$. Such $n$ and $m$ are called the \emph{length} of $a$ and the \emph{length} of $x$, which we denote $|a|$ and $|x|$, respectively.

\noindent We say that $b\in {\cal AR}$ is an end-extension of $a\in {\cal AR}$ and write $a\sqsubseteq b$, if $(\exists B\in \mathcal{R}b=r_n(B))\Rightarrow \exists m\leq |b|$ $(a=r_m(B))$. In an analogous way we define the relation $\sqsubseteq$ on ${\cal AS}$.
  
\begin{itemize}
\item[(A.4)] \textbf{Finitization:} There are relations $\leq_{fin}$ and $\leq^0_{fin}$ on ${\cal AS}\times{\cal AS}$ and ${\cal AR}\times{\cal AS}$, respectively, such that:
\begin{itemize}
\item[(1)] $\{a\colon a\leq^0_{fin}x\}$ and $\{y\colon y\leq_{fin}x\}$ are finite for all $x\in {\cal AS}$.
\item[(2)] $X\leq Y$ iff $\forall n$ $\exists m$ $s_n(X)\leq_{fin}s_m(Y)$.
\item[(3)] $A\leq^0 X$ iff $\forall n$ $\exists m$ $r_n(A)\leq^0_{fin}s_m(X)$.
\item[(4)] $\forall a\in {\cal AR}$ $\forall x,y\in {\cal AS}$ $[a\leq^0_{fin}x\leq_{fin}y\Rightarrow (a\leq^0_{fin}y)]$.
\item[(5)] $\forall a,b\in {\cal AR}$ $\forall x\in {\cal AS}$ $[a\sqsubseteq b$ and $b\leq^0_{fin}x\Rightarrow \exists y\sqsubseteq x$ $(a\leq^0_{fin}y)]$.
\end{itemize}
\end{itemize}

\noindent We deal with the \emph{basic sets} 
$$[a,Y]=\{A\in {\cal R}\colon A\leq^0 Y {\mbox { and }} \exists n \, (r_n(A)=a)\}$$
$$[x,Y]=\{X\in {\cal S}\colon X\leq Y {\mbox { and }} \exists n \, (s_n(X)=x)\}$$
for $a\in {\cal AR}$, $x\in {\cal AS}$ and $Y\in {\cal S}$. Notation: 
$$[n,Y]=[s_n(Y),Y]$$
Also, we define the \emph{depth of } $a\in {\cal AR}$ in $Y\in {\cal S}$ by
$$depth_Y(a)=\left\{
\begin{array}{ll}
min\{k\colon a\leq_{fin}^0 s_k(Y)\}, & {\mbox{if  }} \, \exists k\, (a\leq_{fin}^0 s_k(Y)) \\
-1 , & {\mbox{ otherwise}}
\end{array}\right.$$ 
\noindent The next result is immediate.

\begin{lem}
If $a\sqsubseteq b$ then $depth_Y(a)\leq depth_Y(b)$.  $\blacksquare$
\end{lem}

\noindent Now we state the last two axioms:

\begin{itemize}
\item[(A.5)] \textbf{Amalgamation:} $\forall a\in {\cal AR}$, $\forall \, Y\in {\cal S}$, if $depth_Y(a)=d$, then:
\begin{itemize}
\item[(1)] $d\geq 0\Rightarrow \forall X\in [d,Y]$ $([a,X]\neq \emptyset)$.
\item[(2)] Given $X\in {\cal S}$, 
$$(X\leq Y {\mbox { and }} [a,X]\neq \emptyset) \Rightarrow \exists Y'\in[d,Y]\, ([a,Y']\subseteq [a,X])$$
\end{itemize}
\item[(A.6)] \textbf{Pigeon hole principle:} Suppose $a\in {\cal AR}$ has length $l$ and ${\cal O}\subseteq {\cal AR}_{l+1}=r_{l+1}({\cal R})$. Then for every $Y\in {\cal S}$ with $[a,Y]\neq \emptyset$, there exists $X\in[depth_Y(a), Y]$ such that $r_{l+1}([a,X])\subseteq {\cal O}$ or $r_{l+1}([a,X])\subseteq {\cal O}^{\,c}$.
\end{itemize}

\begin{defn}
We say that ${\cal X}\subseteq {\cal R}$ is ${\cal S}$-Ramsey if for every $[a,Y]$ there exists $X\in[depth_Y(a),Y]$ such that $[a,X]\subseteq {\cal X}$ or $[a,X]\subseteq {\cal X}^c$. If for every $[a,Y]\neq \emptyset$ there exists $X\in[depth_Y(a),Y]$ such that $[a,X]\subseteq {\cal X}^c$, we say that ${\cal X}$ is ${\cal S}$-Ramsey null.
\end{defn}

\begin{defn}
We say that ${\cal X}\subseteq {\cal R}$ is ${\cal S}$-Baire if for every $[a,Y]\neq \emptyset$ there exists a nonempty $[b,X]\subseteq[a,Y]$ such that $[b,X]\subseteq {\cal X}$ or $[b,X]\subseteq {\cal X}^c$. If for every $[a,Y]\neq \emptyset$ there exists a nonempty $[b,X]\subseteq[a,Y]$ such that $[b,X]\subseteq {\cal X}^c$, we say that ${\cal X}$ is ${\cal S}$-meager.  
\end{defn}

\noindent It is clear that every ${\cal S}$-Ramsey set is ${\cal S}$-Baire and every ${\cal S}$-Ramsey null set is ${\cal S}$-meager.

\medskip

\noindent Considering ${\cal AS}$ with the discrete topology and ${\cal AS}^{\mathbb{N}}$ with the completely metri\-zable product topology; we say that ${\cal S}$ is \emph{closed} if it corrresponds to a closed subset of ${\cal AS}^{\mathbb{N}}$ via the identification $X\rightarrow (s_n(X))_{n\in \mathbb{N}}$. 

\begin{defn}
We say that $({\cal R},{\cal S},\leq,\leq^0,r,s)$ is a \emph{Ramsey space} if every ${\cal S}$-Baire subset of ${\cal R}$ is ${\cal S}$-Ramsey and every ${\cal S}$-meager subset of ${\cal R}$ is ${\cal S}$-Ramsey null.
\end{defn}
\begin{thm}[Abstract Ramsey theorem]
Suppose $({\cal R},{\cal S},\leq,\leq^0,r,s)$ sa\-tisfies (A.1)$\dots$(A.6) and ${\cal S}$ is closed. Then  $\blacksquare$
\end{thm} 

\noindent \textbf{Example: The Hales-Jewett space}

\noindent Fix a countable alphabet $L=\cup_{n\in \mathbb{N}}L_n$ with $L_n\subseteq L_{n+1}$ and $L_n$ finite for all $n$; fix $v\notin L$ a "\emph{variable}" and denote $W_L$ and $W_{Lv}$ the semigroups of words over $L$ and of variable words over $L$, respectively. Given $X=(x_n)_{n\in\mathbb{N}}\subseteq W_L\cup W_{Lv}$, we say that $X$ is \emph{rapidly increasing} if 
$$|x_n|>\sum_{i=0}^{n-1}|x_i|$$for all $n\in \mathbb{N}$. Put 
$$W_L^{[\infty]}=\{X=(x_n)_{n\in\mathbb{N}}\subseteq W_L\colon X{\mbox{ is rapidly increasing }}\}$$ 
$$W_{Lv}^{[\infty]}=\{X=(x_n)_{n\in\mathbb{N}}\subseteq W_{Lv}\colon X{\mbox{ is rapidly increasing }}\}$$
By restricting to finite sequences with
$$r_n\colon W_L^{[\infty]}\to W_L^{[n]}\qquad s_n\colon W_{Lv}^{[\infty]}\to W_{Lv}^{[n]}$$we have rapidly increasing finite sequences of words or variable words. The \emph{combinatorial subspaces} are defined for every $X\in W_{Lv}^{[\infty]}$ by
$$[X]_L=\{x_n[\lambda_0]^{\smallfrown} \cdots ^{\smallfrown} x_{n_k}[\lambda_k]\in W_L\colon n_o<\cdots<n_k ,\, \, \lambda_i \in L_{n_i}\}$$
$$[X]_{Lv}=\{x_n[\lambda_0]^{\smallfrown} \cdots ^{\smallfrown} x_{n_k}[\lambda_k]\in W_{Lv}\colon n_o<\cdots<n_k ,\, \, \lambda_i \in L_{n_i}\cup \{v\}\}$$where "$^{\smallfrown}$" denotes concatenation of words and $x[\lambda]$ is the evaluation of the variable word $x$ on the letter $\lambda$.

\noindent For $w\in [X]_L\cup [X]_{Lv}$ we call \emph{support of} $w$ \emph{in} $X$ the unique set $supp_X(w)=\{n_0<n_1<\cdots<n_k\}$ such that $w=x_n[\lambda_0]^{\smallfrown} \cdots ^{\smallfrown} x_{n_k}[\lambda_k]$ as in the definition of the combinatorial subspaces $[X]_L$ and $[X]_{Lv}$. We say that $Y=(y_n)_{n\in \mathbb{N}}\in W_{Lv}^{[\infty]}$ is a \emph{block subsequence} of $X=(x_n)_{n\in \mathbb{N}}\in W_{Lv}^{[\infty]}$ if $y_n\in [X]_{Lv}$ $\forall n$ and
$$max(supp_X(y_n)) < min(supp_X(y_m))$$whenever $n<m$, and write $Y\leq X$. We define the relation $\leq^0$ on $W_L^{[\infty]}\times W_{Lv}^{[\infty]}$ in the natural way. 
Then, if $({\cal R},{\cal S},\leq,\leq^0,r,s)=(W_L^{[\infty]},W_{Lv}^{[\infty]},\leq,\leq^0,r,s)$ is as before, where $r$, $s$ are the restrictions
$$r_n(X)=(x_0,x_1,\dots,x_{n-1})\qquad s_n(Y)=(y_0,y_1,\dots,y_{n-1})$$we have (A.1)$\dots$(A.6), particularly, (A.6) is the well known result:

\begin{thm}
For every finite coloring of $W_L\cup W_{Lv}$ and every $Y\in W_{Lv}^{[\infty]}$ there exists $X\leq Y$ in $W_{Lv}^{[\infty]}$ such that $[X]_L$ and $[X]_{Lv}$ are monochromatic. $\blacksquare$
\end{thm}   

And as a particular case of theorem 1, we have (see \cite{HaJ})

\begin{thm}[Hales--Jewett]
The field of $W_{Lv}^{[\infty]}$-Ramsey subsets of $W_L^{[\infty]}$ is closed under the Souslin operation and it coincides with the field of $W_{Lv}^{[\infty]}$-Baire subsets of $W_L^{[\infty]}$. Moreover, the ideals of $W_{Lv}^{[\infty]}$-Ramsey null subsets of $W_L^{[\infty]}$ and $W_{Lv}^{[\infty]}$-meager subsets of $W_L^{[\infty]}$ are $\sigma$-ideals and they also coincide. $\blacksquare$
\end{thm}

\section{The parametrization}\label{paramet}

\noindent We will denote the family of perfect subsets of $2^{\infty}$ by $\mathbb{P}$ and define

$${\cal AR}[X]= \{ b\in {\cal AR}\colon [ b,X ]\neq \emptyset \}$$also we'll use this notation
$$M\in \mathbb{P}\upharpoonright Q \Leftrightarrow \, (M\in \mathbb{P})\, \wedge \, (M\subseteq Q)$$

\vspace{0.45cm}

\noindent From now on we assume that $({\cal R}, {\cal S}, \leq, \leq^0, r, s )$ is an Ramsey space; that is, (A.1)$\dots$(A.6) hold and ${\cal S}$ is closed. The following are the abstract versions of \emph{perfectly-Ramsey} sets and the $\mathbb{P}\times Exp({\cal R})$-\emph{Baire property} as defined in \cite{mij}.

\vspace{0.45cm}

\begin{defn} $\Lambda\subseteq 2^{\infty}\times {\cal R}$ is \emph{perfectly ${\cal S}$-Ramsey} if for every $Q\in \mathbb{P}$ and $[a,Y]\neq \emptyset$, there exist $M\in \mathbb{P}\upharpoonright Q$ and $X\in [depth_Y(a),Y]$ with $[a,X]\neq \emptyset$ such that $M\times [a,X]\subseteq \Lambda$ or $M\times [a,X]\subseteq \Lambda^c$. If for every $Q\in \mathbb{P}$ and $[a,Y]\neq \emptyset$, there exist $M\in \mathbb{P}\upharpoonright Q$ and $X\in [depth_Y(a),Y]$ with $[a,X]\neq \emptyset$ such that $M\times [a,X]\subseteq \Lambda^c$; we say that  $\Lambda$ is \emph{perfectly ${\cal S}$-Ramsey null}.
\end{defn}

\begin{defn} $\Lambda\subseteq 2^{\infty}\times {\cal R}$ is \emph{perfectly} ${\cal S}$-\emph{Baire} if for every $Q\in \mathbb{P}$ and $[a,Y]\neq \emptyset$, there exist $M\in \mathbb{P}\upharpoonright Q$ and $[b,X]\subseteq[a,Y]$ such that $M\times [b,X]\subseteq \Lambda$ or $M\times [b,X]\subseteq \Lambda^c$. If for every $Q\in \mathbb{P}$ and $[a,Y]\neq \emptyset$, there exist $M\in \mathbb{P}\upharpoonright Q$ and $[b,X]\subseteq[a,Y]$ such that $M\times [b,X]\subseteq \Lambda^c$; we say that $\Lambda$ is \emph{perfectly} ${\cal S}$-\emph{meager}. 
\end{defn}

\noindent Now, the natural extension of combinatorial forcing will be given. From now on fix ${\cal F}\subseteq 2^{<\infty}\times{\cal AR}$ and $\Lambda\subseteq 2^{\infty}\times{\cal R}$.

\vspace{0.45cm}

\noindent \textbf{Combinatorial forcing 1} Given $Q\in \mathbb{P}$, $Y\in {\cal S}$ and $(u,a)\in 2^{<\infty}\times{\cal AR}[Y]$; we say that $(Q,Y)$ \emph{accepts} $(u,a)$ if for every $x\in Q(u)$ and for every $B\in [a,Y]$ there exist integers $k$ and $m$ such that $(x_{|k},r_m(B))\in {\cal F}$. 

\vspace{0.45cm}

\noindent \textbf{Combinatorial forcing 2} Given $Q\in \mathbb{P}$, $Y\in {\cal S}$ and $(u,a)\in 2^{<\infty}\times{\cal AR}[Y]$; we say that $(Q,Y)$ \emph{accepts} $(u,a)$ if $Q(u)\times [a,Y]\subseteq \Lambda$. 

\vspace{0.45cm}

\noindent For both combinatorial forcings we say that $(Q,Y)$ \emph{rejects} $(u,a)$ if for every $M\in \mathbb{P}\upharpoonright Q(u)$ and for every $X\leq Y$ compatible with $a$; $(M,X)$ does not accept $(u,a)$. Also, we say that $(Q,Y)$ \emph{decides} $(u,a)$ if it accepts or rejects it.

\vspace{0.45cm}

\noindent The following lemmas hold for both combinatorial forcings.

\begin{lem}
\begin{itemize}
\item[a)] If $(Q,Y)$ accepts (rejects) $(u,a)$ then $(M,X)$ also accepts (rejects) $(u,a)$ for every $M\in \mathbb{P}\upharpoonright Q(u)$ and for every $X\leq Y$ compatible with $a$.
\item[b)] If $(Q,Y)$ accepts (rejects) $(u,a)$ then $(Q,X)$ also accepts (rejects) $(u,a)$ for every $X\leq Y$ compatible with $a$.
\item[c)] For all $(u,a)$ and $(Q,Y)$ with $[a,Y]\neq\emptyset$, there exist $M\in \mathbb{P}\upharpoonright Q$ and $X\leq Y$ compatible with $a$, such that $(M,X)$ decides $(u,a)$. 
\item[d)] If $(Q,Y)$ accepts $(u,a)$ then $(Q,Y)$ accepts $(u,b)$ for every $b\in r_{|a|+1}([a,Y])$.
\item[e)] If $(Q,Y)$ rejects $(u,a)$ then there exists $X\in[depth_Y(a),Y]$ such that $(Q,Y)$ does not accept $(u,b)$ for every $b\in r_{|a|+1}([a,X])$. 
\item[f)] $(Q,Y)$ accepts (rejects) $(u,a)$ iff $(Q,Y)$ accepts (rejects) $(v,a)$ for every $v\in 2^{<\infty}$ such that $u\sqsubseteq v$.
\end{itemize}
\end{lem}
\noindent \textbf{Proof:} (a) and (b) follow from the inclusion: $M(u)\times [a,X]\subseteq Q(u)\times [a,Y]$ if $X\leq Y$ and $M\subseteq Q(u)$.

\vspace{0.3cm}

\noindent(c) Suppose that we have $(Q,Y)$ such that for every $M\in \mathbb{P}\upharpoonright Q$ and every $X\leq Y$ compatible with $a$, $(M,X)$ does not decide $(u,a)$. Then $(M,X)$ does not accept $(u,a)$ if $M\in \mathbb{P}\upharpoonright Q(u)$; i.e. $(Q,Y)$ rejects $(u,a)$. 

\vspace{0.3cm}

\noindent (d) Follows from: $a\sqsubseteq b$ and $[a,Y]\subseteq [b,Y]$, if $b\in r_{|a|+1}([a,X])$.

\vspace{0.3cm}

\noindent (e) Suppose $(Q,Y)$ rejects $(u,a)$ and define $\phi\colon {\cal AR}_{|a|+1}\to 2$ by $\phi(b)=1$ if $(Q,Y)$ accepts $(u,b)$. By (A.6) there exist $X\in[depth_Y(a),Y]$ such that $\phi$ is constant in $r_{|a|+1}([a,X])$. If $\phi(r_{|a|+1}([a,X]))=1$ then $(Q,X)$ accepts $(u,a)$, which contradicts $(Q,Y)$ rejects $(u,a)$ (by part (b)). The result follows.

\vspace{0.3cm}

\noindent(f) ($\Leftarrow$)Obvious.

\vspace{0.3cm}

\noindent ($\Rightarrow$) Follows from the inclusion: $Q(v)\subseteq Q(u)$ if $u\sqsubseteq v$. \quad $\blacksquare$

\vspace{0.45cm}

\noindent We say that a sequence $([n_k,Y_k])_{k\in \mathbb{N}}$ is a \emph{fusion sequence} if:
\begin{enumerate}
\item $(n_k)_{k\in \mathbb{N}}$ is nondecreasing and converges to $\infty$.
\item $X_{k+1}\in [n_k,X_{k}]$ for all $k$.
\end{enumerate}

\noindent Note that since ${\cal S}$ is closed, for every fusion sequence $([n_k,Y_k])_k\in \mathbb{N}$ there exist a unique $Y \in {\cal S}$ such that $s_{n_k}(Y)=s_{n_k}(X_k)$ and $Y \in [n_k,X_k]$ for all $k$. $Y$ is called the \emph{fusion} of the sequence and is denoted $lim_kX_k$.

\begin{lem}
Given $P\in \mathbb{P}$, $Y\in {\cal S}$ and $N\geq0$; there exist $Q\in \mathbb{P}\upharpoonright P$ and $X\leq Y$ such that $(Q,X)$ decides every $(u,a)\in 2^{<\infty}\times {\cal AR}[X]$ with $N\leq depth_X(a)\leq |u|$.
\end{lem}

\noindent \textbf{Proof:} We build sequences $(Q_k)_k$ and $(Y_k)_k$ such that:
\begin{enumerate}
\item $Q_0=P$, $Y_0=Y$.
\item $n_k=N+k$.
\item $(Q_{k+1},Y_{k+1})$ decides every $(u,b)\in 2^{n_k}\times {\cal AR}[Y_k]$ with $depth_{Y_k}(b)=n_k$.
\end{enumerate}

\noindent Suppose we have defined $(Q_k,Y_k)$. List $\{ b_0, \dots, b_r\}=\{b\in {\cal AR}[Y_k]\colon depth_{Y_k}(b)=n_k\}$ and $\{u_0,\dots,u_{2^{n_k}-1}\}=2^{n_k}$. By lemma 1(c) there exist $Q_k^{0,0}\in \mathbb{P}\upharpoonright Q_k(u_0)$ and $Y_k^{0,0}\in[n_k,Y_k]$ compatible with $b_0$ such that $(Q_k^{0,0},Y_k^{0,0})$ decides $(u_0,b_0)$. In this way we can obtain $(Q_k^{i,j},Y_k^{i,j})$ for every $(i,j)\in \{0,\dots,2^{n_k}-1\}\times\{0,\dots,r\}$, which decides $(u_i,b_j)$ and such that $Q_k^{i,j+1}\in \mathbb{P}\upharpoonright Q_k^{i,j}(u_i)$, $Y_k^{i,j+1}\leq Y_k^{i,j}$ is compatible with $b_{j+1}$, $Q_k^{i+1,0}\in \mathbb{P}\upharpoonright Q_k(u_{i+1})$ and $Y_k^{i+1,0}\leq Y_k^{i,r}$.  

\noindent Define $$Q_{k+1}=\bigcup_{i=0}^{2^{n_k-1}}Q_k^{i,r}\qquad , \qquad Y_{k+1}=Y_k^{2^{n_k-1},r}$$

\noindent Then, given $(u,b)\in 2^{n_k}\times {\cal AR}[Y_{k+1}]$ with  $depth_{Y_{k+1}}(b)=n_k=depth_{Y_k}(b)$, there exist $(i,j)\in \{0,\dots,2^{n_k}-1\}\times \{0,\dots,r\}$ such that $u=u_i$ and $b=b_j$. So $(Q_k^{i,j},Y_k^{i,j})$ decides $(u,b)$ and, since
$$Q_{k+1}(u_i)=Q_k^{i,r}\subseteq Q_k^{i,j}(u_i)\subseteq Q_k^{i,j} {\mbox { and } }Y_{k+1}\leq Y_k$$
we have $(Q_{k+1},Y_{k+1})$ decides $(u,b)$ (by lemma 1(a)) We claim that $Q=\cap_k Q_k$ and $X=lim_k Y_k$ are as required: given $(u,a)\in 2^{<\infty}\times {\cal AR}[X]$ with $N\leq depth_X(a)\leq |u|$, we have $depth_X(a)=n_k=depth_{Y_k}(a)$ for some $k$. Then, if $|u|=n_k$, $(Q_{k+1},Y_{k+1})$ from the construction of $X$ decides $(u,a)$ and hence $(Q,X)$ decides $(u,a)$. If $|u|>n_k$ $(Q,X)$ decides $(u,a)$ by lemma 1(f).  $\quad \blacksquare$   

\begin{lem}
Given $P\in \mathbb{P}$, $Y\in {\cal S}$, $(u,a)\in 2^{<\infty}\times {\cal AR}[Y]$ with $depth_Y(a)\leq |u|$ and $(Q,X)$ as in lemma 2 with $N=depth_Y(a)$; if $(Q,X)$ rejects $(u,a)$ then there exist $Z\leq X$ such that $(Q,Z)$ rejects $(v,b)$ if $u\sqsubseteq v$, $a\sqsubseteq b$ and $depth_Z(b)\leq |v|$.
\end{lem}

\noindent \textbf{Proof:} Let's build a fusion sequence $([n_k,Z_k])_k$, with $n_k=|u|+k$. Let $Z_0=X$. Then $(Q,Z_0)$ rejects $(u,a)$ (and by lemma 1(f) it rejects $(v,a)$ if $u\sqsubseteq v$). Suppose we have $(Q,Z_k)$ which rejects every $(v,b)$ with $v\in 2^{n_k}$ extending $u$, $a\sqsubseteq b$ and $depth_{Z_k}(b)\leq n_k$. List $\{b_0,\dots,b_r\}=\{b\in {\cal AR}[Z_k]\colon a\sqsubseteq b {\mbox { and }} depth_{Z_k}(b)\leq n_k\}$ and $\{u_0,\dots,u_s\}$ the set of all $v\in 2^{n_k+1}$ extending $u$. By lemma 1(f) $(Q,Z_k)$ rejects $(u_i,b_j)$, for every $(i,j)\in \{0,\dots,s\}\times\{0,\dots,r\}$. Use lemma 1(e)  to find $Z_k^{0,0}\in[n_k,Z_k]$ such that $(Q,Z_k^{0,0})$ rejects $(u_0,b)$ if $b\in r_{|b_0|+1}([b_0,Z_k^{0,0}])$. In this way, for every $(i,j)\in \{0,\dots,s\}\times\{0,\dots,r\}$, we can find $Z_k^{i,j}\in[n_k,Z_k]$ such that $Z_k^{i,j+1}\in [n_k,Z_k^{i,j}]$, $Z_k^{i+1,0}\in [n_k,Z_k^{i,r}]$ and $(Q,Z_k^{i,j})$ rejects $(u_i,b)$ if $b\in r_{|b_j|+1}([b_j,Z_k^{i,j}])$. Define $Z_{k+1}=Z_k^{s,r}$. Note that if $(v,b)\in 2^{<\infty}\times {\cal AR}[Z_{k+1}]$, $a\sqsubseteq b$, $u\sqsubseteq v$ and $depth_{Z_{k+1}}(b)=n_k+1$ then $v=u_i$ for some $i\in \{0,\dots,s\}$ and $b=r_{|b|}(A)$, $a=r_{|a|}(A)$ for some $A\leq^0 Z_{k+1}$; by (A.4)(5) there exist $m\leq n_k$ such that $b'=r_{|b|-1}(A)\leq^0_{fin}s_m(Z_{k+1})$, so $depth_{Z_{k+1}}(b')\leq n_k$, i.e. $b'=b_j$ for some $j\in \{0,\dots,r\}$. Then $b\in r_{|b_j|+1}([b_j,Z_k^{i,j}])$. Hence, by lemma 1(f), $(Q,Z_{k+1})$ rejects $(v,b)$. Then $Z=lim_k Z_k$ is as required: given $(v,b)$ with $u\sqsubseteq v$, $a\sqsubseteq b$ and $depth_Z(b)\leq |v|$ then $depth_Z(b)=depth_Y(a)+k\leq n_k$ for some $k$ and $b\in r_{|b_j|+1}([b_j,Z_k^{i,j}])$ for some $j\in \{0,\dots,r\}$ from the construction of $Z$ (again, by (A.4)(5)). So $(Q,Z_k)$ (from the construction of $Z$) rejects $(v,b)$ and, by lemma 1(a), $(Q,Z)$ also does it. \quad $\blacksquare$ 

\vspace{0.45cm}

\noindent The following theorem is an extension of theorem 3 \cite{mij} and its proof is analogous. 

\begin{thm}\label{perf-galvin}
For every ${\cal F}\subseteq 2^{<\infty}\times {\cal AR}$, $P\in \mathbb{P}$, $Y\in {\cal S}$ and $(u,a)\in 2^{<\infty}\times {\cal AR}$ there exist $Q\in \mathbb{P}\upharpoonright P$ and $X\leq Y$ such that one of the following holds:
\begin{enumerate}
\item For every $x\in Q$ and $A\in [a,X]$ there exist integers $k$, $m>0$ such that $(x_{|k}, r_m(A))\in {\cal F}$.
\item $(T_Q\times {\cal AR}[X])\cap {\cal F}=\emptyset$. 
\end{enumerate} 
\end{thm} 

\noindent \textbf{Proof:} Whitout loss of generality, we can assume $(u,a)=(\left\langle  \right\rangle,\emptyset)$. Consider combinatorial forcing 1. Let $(Q,X)$ as in lemma 3 ($N=0$). If
$(Q,X)$ accepts $(\left\langle  \right\rangle,\emptyset)$, part (1) holds. If $(Q,X)$ rejects $(\left\langle  \right\rangle,\emptyset)$, use lemma 4 to obtain $Z\leq X$ such that $(Q,X)$ rejects $(u,a)\in2^{<\infty}\times{\cal AR}[Z]$ if $depth_Z(b)\leq|u|$. If $(t,b)\in(T_Q\times {\cal AR}[Z])\cap {\cal F}$, find $u_t\in2^{<\infty}$ such that $Q(u_t)\subseteq Q\cap [t]$. Thus, $(Q,Z)$ accepts $(u,b)$. In fact: for $x\in Q(u_t)$ and $B\in [b,Z]$ we have $(x_{|k}, r_m(A))=(t,b)\in {\cal F}$ if $k=|t|$ and $m=|b|$. By lemma 2(f), $(Q,Z)$ accepts $(v,b)$ if $u_t\sqsubseteq v$ and $depth_Z(b)\leq|v|$. But this is a contradiction with the choice of $Z$. Hence, $(T_Q\times {\cal AR}[X])\cap {\cal F}=\emptyset$. $\blacksquare$

The next theorem is our main result and its proof is analogous to theorem 3 \cite{mij}. 

\begin{thm}\label{main}
For $\Lambda\subseteq 2^{\infty}\times {\cal R}$ we have:
\begin{enumerate}
\item $\Lambda$ is perfectly ${\cal S}$-Ramsey iff it is perfectly ${\cal S}$-Baire.
\item $\Lambda$ is perfectly ${\cal S}$-Ramsey null iff it is perfectly ${\cal S}$-meager. 
\end{enumerate} 
\end{thm} 

\noindent \textbf{Proof:} (1) We only have to prove the implication from right to left. Suppose that $\Lambda\subseteq 2^{\infty}\times {\cal R}$ is perfectly ${\cal S}$-Baire. Again, whitout loss of generality, we can lead whith a given $Q\times[\emptyset,Y]$. Using combinatorial forcing and lemma 3, we have the following:

\begin{clm}
Given $\hat{\Lambda}\subseteq 2^{\infty}\times {\cal R}$, $P\in \mathbb{P}$ and $Y\in{\cal S}$, there exists $Q\in \mathbb{P}\upharpoonright P$ and $X\leq Y$ such that for each $(u,b)\in 2^{<\infty}\times {\cal AR}[X]$ with $depth_X(b)\leq |u|$ one of the following holds:
\begin{itemize}
\item[i)] $Q(u)\times[b,X]\subseteq \hat{\Lambda}$
\item[ii)] $R\times[b,Z]\not\subseteq \hat{\Lambda}$ for every $R\subseteq Q(u)$ and every $Z\leq X$ compatible with $b$. 
\end{itemize}
\end{clm}

By applying the claim to $\Lambda$, $P$ and $Y$, we find $Q_1\in \mathbb{P}\upharpoonright P$ and $X_1\leq Y$ such that for each $(u,b)\in 2^{<\infty}\times {\cal AR}[X_1]$ with $depth_{X_1}(b)\leq |u|$ one of the following holds:
\begin{itemize}
\item $Q_1(u)\times[b,X_1]\subseteq \Lambda$ or
\item $R\times[b,Z]\not\subseteq \Lambda$ for every $R\subseteq Q_1(u)$ and every $Z\leq X_1$ compatible with $b$. 
\end{itemize}

For each $t\in T_{Q_1}$, choose $u_1^t\in 2^{<\infty}$ with $u_1^t(Q_1)\sqsubseteq t$. If we define the family 
$${\cal F}_1=\{(t,b)\in T_{Q_1}\times {\cal AR}[X_1]\colon Q_1(u_1^t)\times [b,X_1]\subseteq\Lambda\}$$then we find $S_1\subseteq Q_1$and $Z_1\leq X_1$ as in theorem 4. If (1) of theorem \ref{perf-galvin} holds, we are done. If part (2) holds, apply the claim to $\Lambda^c$, $S_1$ and $Z_1$ to find $Q_2\in \mathbb{P}\upharpoonright P$ and $X_2\leq Y$ such that for each $(u,b)\in 2^{<\infty}\times {\cal AR}[X_2]$ with $depth_{X_2}(b)\leq |u|$ one of the following holds:
\begin{itemize}
\item $Q_2(u)\times[b,X_2]\subseteq \Lambda^c$ or
\item $R\times[b,Z]\not\subseteq \Lambda^c$ for every $R\subseteq Q_2(u)$ and every $Z\leq X_2$ compatible with $b$. 
\end{itemize}  

Again, for each $t\in T_{Q_2}$, choose $u_2^t\in 2^{<\infty}$ with $u_2^t(Q_2)\sqsubseteq t$; define the family 
$${\cal F}_2=\{(t,b)\in T_{Q_2}\times {\cal AR}[X_2]\colon Q_2(u_2^t)\times [b,X_1]\subseteq\Lambda\}$$and find $S_2\subseteq Q_2$ and $Z_2\leq X_2$ as in theorem 4. If (1) holds, we are done and part (2) is not possible since $\Lambda$ is perfectly ${\cal S}$-Baire (see \cite{mij}). This proves (1). To see part (2), notice that, as before, we only have to prove the implication from right to left, which follows from part (1) f $\Lambda$ is perfectly ${\cal S}$-meager.  $\blacksquare$

\begin{coro}[Parametrized infinite dimensional Hales-Jewett theorem]\label{paraHJ}
For $\Lambda\subseteq 2^{\infty}\times W_L^{[\infty]}$ we have:
	\begin{enumerate}
		\item $\Lambda$ is perfectly Ramsey iff it has the $\mathbb{P}\times W_{Lv}^{[\infty]}$-Baire property.
		\item $\Lambda$ is perfectly Ramsey null iff it is $\mathbb{P}\times W_{Lv}^{[\infty]}$-meager . $\blacksquare$
	\end{enumerate}

\end{coro}

\noindent Making $\mathcal{R} = \mathcal{S}$ in $({\cal R}, {\cal S}, \leq, \leq^0, r, s )$, we obtain the following:

\begin{coro}[Mijares]
If $({\cal R},\leq,(p_n)_{n\in\mathbb{N}})$ is a topological Ramsey space then:
\begin{enumerate}
\item $\Lambda\subseteq {\cal R}$ is perfectly Ramsey iff has the $\mathbb{P}\times Exp({\cal R})$-Baire property.
\item $\Lambda\subseteq {\cal R}$ is perfectly Ramsey null iff is $\mathbb{P}\times Exp({\cal R})$-meager. $\blacksquare$
\end{enumerate} 
\end{coro}

\begin{coro}[Pawlikowski]
For $\Delta\subseteq 2^{\infty}\times \mathbb{N}^{[\infty]}$ we have:
	\begin{enumerate}
		\item $\Lambda$ is perfectly Ramsey iff it has the $\mathbb{P}\times Exp(\mathbb{N}^{[\infty]})$-Baire property.
		\item $\Lambda$ is perfectly Ramsey null iff it is $\mathbb{P}\times Exp(\mathbb{N}^{[\infty]})$-meager . $\blacksquare$
	\end{enumerate}
\end{coro}

 \noindent Now we will proof that the family of perfectly ${\cal S}$-Ramsey and perfectly ${\cal S}$-Ramsey null subsets of $2^{\infty}\times {\cal R}$ are closed under the Souslin operation. Recall that the result of applying the Souslin operation to a given $(\Lambda_a)_{a\in\mathcal{AR}}$ is
$$\bigcup_{A\in\mathcal{R}}\bigcap_{n\in\mathbb{N}}\Lambda_{r_n(A)}$$

\begin{prop}
The perfeclty ${\cal S}$-Ramsey null subsets of $2^{\infty}\times{\cal R}$ form a $\sigma$-ideal. 
\end{prop}

\noindent \textbf{Proof:} This proof is also analogous to its corresponding version in \cite{mij} (lemma 4). So we just expose the mean ideas. Given an increasing sequence of perfectly ${\cal S}$-Ramsey null subsets of $2^{\infty}\times{\cal R}$ and $P\times[\emptyset,Y]$, we proceed as in lemma 3 to build fusion sequences $(Q_n)_n$ and $[n+1,X_n]$ such that $$Q_n\times[b,X_n]\cap \Lambda_n=\emptyset$$ for every $n\in \mathbb{N}$ and $b\in {\cal AR}[X_n]$ with $depth_{X_n}(b)\leq n$. Thus, if $Q=\cap_n Q_n$ and $X=lim_n X_n$, we have $Q\times [\emptyset,X]\cap \bigcup_n \Lambda_n=\emptyset$. $\blacksquare$

\medskip

\noindent Recall that given a set $X$, two subsets $A,B$ of $X$ are "\emph{compatibles}" with respect to a family ${\cal F}$ of subsets $X$ if there exists $C\in {\cal F}$ such that $C\subseteq A\cap B$. And ${\cal F}$ is \emph{M-like} if for ${\cal G}\subseteq {\cal F}$ with $|{\cal G}|<|{\cal F}|$, every member of ${\cal F}$ which is not compatible with any member of ${\cal G}$ is compatible with $X\setminus \bigcup {\cal G}$. A $\sigma$-algebra ${\cal A}$ of subsets of $X$ together with a $\sigma$-ideal ${\cal A}_0\subseteq{\cal A}$ is a \emph{Marczewski pair} if for every $A\subseteq X$ there exists $\Phi(A)\in {\cal A}$ such that $A\subseteq \Phi(A)$ and for every $B\subseteq \Phi(A)\setminus A$, $B\in{\cal A}\rightarrow B\in {\cal A}_0$. The following is a well known fact:

\begin{thm}[Marczewski]\label{marcz}
Every $\sigma$-algebra of sets which together with a $\sigma$-ideal is a Marczeswki pair, is closed under the Souslin operation.
$\quad \blacksquare$
\end{thm} 

Let's denote ${\cal E}({\cal S})=\{[n,Y]\colon n\in \mathbb{N}$, $Y\in {\cal S}\}$.

\begin{prop}
If $|{\cal S}|=2^{\aleph_0}$, then the family ${\cal E}({\cal S})$ is $M$-like.
\end{prop}

\noindent \textbf{Proof:} Consider ${\cal B}\subseteq {\cal E}({\cal S})$ with $|{\cal B}|<|{\cal E}({\cal S})|=2^{\aleph_0}$ and suppose that $[a,Y]$ is not compatible with any member of ${\cal B}$, i. e. for every $B\in {\cal B}$, $B\cap [a,Y]$ does not contain any member of ${\cal E}({\cal S})$. We claim that $(Q,Y)$ is compatible with ${\cal R}\smallsetminus \bigcup {\cal B}$. In fact:

\noindent Since $|{\cal B}|<2^{\aleph_0}$, $\bigcup {\cal B}$ is ${\cal S}$-Baire  (it is ${\cal S}$-Ramsey). So, there exist $[b,X]\subseteq [a,Y]$ such that:
\begin{enumerate}
\item $[b,X]\subseteq \bigcup {\cal B}$ or
\item $[b,X]\subseteq {\cal R}\smallsetminus \bigcup {\cal B}$
\end{enumerate}

\noindent(1) is not possible because $[a,Y]$ is not compatible with any member of ${\cal B}$. And (2) says that $[a,Y]$ is compatible with ${\cal R}\smallsetminus \bigcup {\cal B}$ \quad{$\blacksquare$} 

As consequences of the previous proposition and theorem \ref{marcz}, the following facts hold.

\begin{coro}\label{souslin}
If $|{\cal S}|=2^{\aleph_0}$, then the family of perfectly ${\cal S}$-Ramsey subsets of $2^{\infty}\times {\cal R}$ is closed under the Souslin operation. $\blacksquare$
\end{coro} 

\begin{coro}
The field of perfectly $W_{Lv}^{[\infty]}$-Ramsey subsets of $2^{\infty}\times W_L^{[\infty]}$ is closed under the Souslin operation. $\blacksquare$
\end{coro}

\noindent

\medskip

\noindent Finally, making $\mathcal{R} = \mathcal{S}$ in $({\cal R}, {\cal S}, \leq, \leq^0, r, s )$, we obtain the following:

\begin{coro}[Mijares]
If $({\cal R},\leq,r)$ satisfies (A.1)$\dots$(A.6), ${\cal R}$
is closed, and $|{\cal R}|=2^{\aleph_0}$ then the family of perfectly Ramsey subsets of $2^{\infty}\times {\cal R}$ is closed under the Souslin operation. $\blacksquare$
\end{coro}

\begin{coro}[Pawlikowski]
The field of perfectly Ramsey subsets of $2^{\infty}\times \mathbb{N}^{[\infty]}$ is closed under the Souslin operation. $\blacksquare$
\end{coro}

\label{end-art}

\begin{thebibliography}{99}
\bibitem{carsimp} Carlson, T. J, Simpson, S. G. {\em Topological Ramsey theory}, in Ne\^{s}et\^{r}il, J., R\"{o}dl, {\em Mathematics of Ramsey Theory}(Eds.), Springer, Berlin, 1990, pp. 172--183.
\bibitem{dip} Di Prisco, C., {\em Partition properties and perfect sets}, Adv. in Math., \textbf{176}(2003), 145--173.
\bibitem{diptodo} Di Prisco, C., Todorcevic, S., {\em Souslin partitions of products of finite sets}, Notas de L\'ogica Matem\'atica Vol. 8, Universidad Nacional del Sur, Bah\'ia Blanca, Argentina, 1993, pp. 119-127.
\bibitem{ellen} Elentuck, E. {\em A new proof that analitic sets are Ramsey}, J. Symbolic Logic, \textbf{39}(1974), 163--165.
\bibitem{farah} Farah, I. {\em Semiselective coideals}, Mathematika., \textbf{45}(1998), 79--103.
\bibitem{galpri} Galvin, F., Prikry, K. {\em Borel sets and Ramsey's theorem}, J. Symbolic Logic, \textbf{38}(1973), 193--198.
\bibitem{HaJ} Hales, A.W. and Jewett, R.I., {\em Regularity and Propositional
Games}, Trans. Amer. Math. \textbf{106} (1963), 222-229.
\bibitem{mij} Mijares, J. {\em Parametrizing the abstract Ellentuck theorem}, Discrete Math., \textbf{307}(2007), 216--225.
\bibitem{miltodo} A. Miller, {\em Infinite combinatorics and definibility}, Ann. Pure Appl. Logic  \textbf{41}(1989), 178--203.
\bibitem{mill} Milliken, K., \emph{Ramsey's theorem with sums or unions}, J. Comb. Theory, ser A \textbf{18}(1975), 276--290.
\bibitem{nash} Nash-Williams, C. St. J. A., {\em On well-quasi-ordering transfinite sequences}, Proc. Cambridge Philo. Soc., \textbf{61}(1965), 33--39.
\bibitem{pawl} J. Pawlikowski, {\em Parametrized Elletuck theorem}, Topology and its applications \textbf{37}(1990), 65--73.
\bibitem{todo}  Todorcevic, S., {\em Introduction to Ramsey spaces},
Princeton University Press, Princeton, New Jersey, 2010.

\end{thebibliography}
\end{document}